\documentclass[9pt]{article}
\usepackage{amsthm,amsfonts,amsmath,amssymb}
\usepackage[all]{xy}
\theoremstyle{definition}
\newtheorem{thm_definition}{Theorem} [section]

\newtheorem{cor_definition}[thm_definition]{Corollary}
\newtheorem{lem_definition}[thm_definition]{Lemma}
\newtheorem{pro_definition}[thm_definition]{Proposition}

\title{\textbf{Intersection of stable and unstable manifolds for invariant Morse functions   }}
\author{Hitoshi Yamanaka}
\date{ }

\begin{document}
\maketitle
\begin{abstract}
We study the structure of the smooth manifold which is defined as the intersection of a stable manifold and an unstable manifold for an invariant Morse-Smale function.

\end{abstract}

\section{Introduction}
The aim of this paper is to investigate invariant Morse functions on compact smooth manifolds with action of compact Lie groups.
  
Let $M$ be a compact $n$-dimensional Riemannian manifold, $\left<\cdot ,\cdot \right>$ its Riemannian metric,  and  $\Phi $ a Morse function on $M$. We denote by $-\nabla \Phi $ the negative gradient vector field of $\Phi $ with respect to the metric $\left<\cdot, \cdot \right>$, and  let $\gamma _{p}(t)$ be the corresponding negative gradient flow which passes through a point $p$ of $M$ at $t=0$. 
For a critical point $p$ of $\Phi $,  the \textbf{unstable manifold} and the \textbf{stable manifold} of $p$ are defined by
\begin{center}
$W^{u}(p)=\bigg\{ x\in M \bigg| \lim\limits_{t\rightarrow -\infty}\gamma _{x}(t)=p \bigg\}$, \\ 

$W^{s}(p)=\bigg\{ x \in M \bigg| \lim\limits_{t\rightarrow \infty}\gamma _{x}(t)=p \bigg\}$
\end{center}
respectively.
Since $\Phi $ is a Morse function, $W^{u}(p)$ and  $W^{s}(p)$ are a smoothly embedded open disks of dimensions $n-\lambda(p),\lambda (p)$ respectively,
where $\lambda (p)$ denotes the Morse index of $p$ (see [BH, Theorem 4.2]).
We say that a Morse function $\Phi $ is \textbf{Morse-Smale} if $W^{u}(p)$ and $W^{s}(q)$ intersect transversally for all critical points $p,q$.
If the Morse function $\Phi $ is Morse-Smale, then $\widetilde{\mathcal{M}}(p,q):=W^{u}(p)\cap W^{s}(q)$ is also a submanifold of $M$ which has dimension $\lambda (p)-\lambda (q)$. 

$\widetilde{\mathcal{M}}(p,q)$ has a natural $\mathbb{R}$-action which is defined by $t\cdot x:=\gamma _{x}(t)$ where $t\in \mathbb{R},x\in \widetilde{\mathcal{M}}(p,q)$.\ The quotient space of $\widetilde{\mathcal{M}}(p,q)$
by the $\mathbb{R}$-action is denoted by $\mathcal{M}(p,q)$. Witten's Morse theory [W] asserts that in some cases, the homology group of $M$ with integral coefficient is
recovered from the structure of $\mathcal{M}(p,q)$'s such that $\lambda (p)-\lambda (q)=1$.  However, there is a Morse function which has no critical points $p,q$ such that $\lambda (p)-\lambda (q)=1$. For example, for a certain Morse function on the partial flag manifold, every unstable manifold is given by the Bruhat cell $BwP/P$. In particular,
every Morse index is even (see [A]).

This phenomenon leads us to the study of the structure of $\widetilde{\mathcal{M}}(p,q)$ for $p,q\in \operatorname{Cr}(\Phi ),\lambda (p)-\lambda (q)=2$.

In this paper, we investigate the structure of $\mathcal{M}(p,q)$ for $p,q\in \operatorname{Cr}(\Phi )$ such that $\lambda (p)-\lambda (q)=2$ under the assumption
that $M$ admits an action of a compact Lie group $G$ and $\Phi $ is $G$-invariant.

Our main theorem is the following.
\begin{thm_definition}
Let $\Phi $ be a $G$-invariant Bott-Morse function on $M$. Let $p,q$ be $G$-fixed points. Assume the following conditions:\\
(1) $M^{G}\subset \operatorname{Cr}(\Phi )$.\\
(2) $\lambda (p)-\lambda (q)=2$.\\
(3) $W^{u}(p)$ and $W^{s}(q)$ intersect transversally.\\ 
Then every connected component of $\widetilde{\mathcal{M}}(p,q)$ is diffeomorphic to $S^{1}\times \mathbb{R}$.$\hspace{6mm}$ $\square$
\end{thm_definition}
We also show that the action of $G$ on $\widetilde{\mathcal{M}}(p,q)$ is given by the rotation of sphere
(see Proposition 3.6 below). 
By these results geometric structure of $\widetilde{\mathcal{M}}(p,q)$ in our setting is similar to the one treated in the GKM theory [GKM].

This paper is organized as follows. In Section 2, we study the critical point set of an invariant Morse function and apply it to an invariant Morse function on a
homogenious space.\ In Section 3,  we prove Theorem 1.1.

\section{Critical points}

Let $G$ be a compact Lie group and $M$ be a compact $G$-manifold. Denote by $M^{G}$ the fixed point set of
the action of $G$ on $M$. We say a smooth function $\Phi :M\longrightarrow \mathbb{R}$ is $G$-\textbf{invariant} if it satisfies
$\Phi (g\cdot p)=\Phi (p)$ for all $g\in G,p\in M$. For a smooth function $\Phi$ on $M$, we denote by $\operatorname{Cr}(\Phi )$ the critical point
set of $\Phi $.
\begin{pro_definition}
Let $G$ be a compact connected Lie group, $M$ be a compact smooth $G$-manifold, and $\Phi:M\longrightarrow \mathbb{R}$ be a $G$-invariant Morse function on $M$.
Assume that there exist only finitely many $G$-fixed points on $M$. Then we have $\operatorname{Cr}(\Phi )=M^{G}$.
$\hspace{95mm}$ $\square$ 
\label{crit} 
\end{pro_definition}
Since $G$ and $M$ are both compact,\ there exists a $G$-invariant metric $\left<\cdot ,\cdot \right>$ on $M$.
Consider the negative gradient flow equation 
\begin{center}
$\displaystyle \gamma (0)=p,\ \ \frac{d}{dt}\gamma (t)=-(\nabla \Phi )_{\gamma (t)}$.
\end{center}
Here,\ we denote by $\nabla \Phi $ the gradient vector field for $\Phi $ with respect to
the $G$-invariant Riemannian metric $\left<\cdot ,\cdot \right>$ on $M$. Let $\gamma_{p}(t)$ be the unique solution of this equation. By the uniqueness of the solution we see easily the following.
\begin{lem_definition}
We have $\gamma_{g\cdot p}(t)=g\cdot \gamma_{p}(t)$ for all $g\in G,p\in M$.$\hspace{23mm}$ $\square$
\label{ggamma}
\end{lem_definition}
\begin{flushleft} \textbf{Proof of Proposition 2.1.}\ Take $p\in \operatorname{Cr}(\Phi )$.\ By Lemma \ref{ggamma}, we have\end{flushleft}
\begin{center}
$\lim\limits_{t\rightarrow -\infty}\gamma_{g\cdot p}(t)=\lim\limits_{t\rightarrow -\infty}g\cdot \gamma_{p}(t)=g\cdot p$.    
\end{center}
This means $g\cdot p$ is also a critical point for $\Phi $, so we have $G\cdot p\subset \operatorname{Cr}(\Phi )$.
However,\ since $M$ is compact,\ $\operatorname{Cr}(\Phi )$ is a finite set.\ Thus by the
connectedness of $G$,\ we have $G\cdot p=\{ p\}$.\ This shows $p\in M^{G}$.\par
Take $p\in M^{G}$.\ By Lemma \ref{ggamma} we have
\begin{center}
$g\cdot \gamma_{p}(t)=\gamma_{g\cdot p}(t)=\gamma_{p}(t)$
\end{center}
for all $g\in G$.\ This means $\{\gamma _{p}(t) |t\in \mathbb{R} \}\subset M^{G}$.
Since $M^{G}$ is a finite set, this implies $\{\gamma _{p}(t) |t\in \mathbb{R} \}=\{p\}$.\ Thus we have $p\in \operatorname{Cr}(\Phi )$.$\hspace{35mm}$ $\square$

\begin{cor_definition}
Let $p_{0}$ be a point of $M$ and $H$ be its stabilizer.\ Assume the following
three conditions:\\
(1) $H$ is connected.\\
(2) $W_{H}:=N_{G}(H)/H$ is a finite group.\\
(3) The Fixed point set of the $H$-action on $M$ is contained in the $G$-orbit of $p_{0}$.
Then, we have
\begin{center}
$\operatorname{Cr}(\Phi )=W_{H}\cdot p_{0}$
\end{center}
for any $H$-invariant Morse function $\Phi:M\longrightarrow  \mathbb{R}$.
\begin{proof}
First, we prove $M^{H}=W_{H}\cdot p_{0}$. The inclusion $M^{H}\supset W_{H}\cdot p_0$ is clear. Take $p\in M^{H}$. Then by
the condition (3),\ it is contained in the $G$-orbit of $p_{0}$. So we can write $p=g\cdot p_{0}$ where $g$ is an element of $G$. Since $p\in M^{H}$,\ we have
$h\cdot (g\cdot p_{0})=g\cdot p_{0}$ for all $h\in H$. So we have $g^{-1}Hg\subset H$.\ Since $g^{-1}Hg$
and $H$ are connected Lie subgroups with the same Lie algebra, the inclusion
implies $g^{-1}Hg=H$. Thus we have $p=g\cdot p_{0}\in W_{H}\cdot p_{0}$, as desired.

In particular, by the condition (2), $M^{H}=W_{H}\cdot p_{0}$ is a finite set. Thus by Proposition \ref{crit},
we have $\operatorname{Cr}(\Phi )=W_{H}\cdot p_{0}$. 
\end{proof}
\label{nonhomo}\end{cor_definition}
As an application to homogeneous spaces,\ we have the following corollaries:
\begin{cor_definition}
Let $G$ be a compact Lie group and $H$ be its connected closed subgroup.
If $N_{G}(H)/H$ is a finite group, we have
\begin{center}
$\operatorname{Cr}(\Phi)=N_{G}(H)/H$
\end{center}
for any $H$-invariant Morse function $\Phi:G/H\longrightarrow \mathbb{R}$. $\hspace{36mm}$ $\square$
\label{homo}
\end{cor_definition}
\begin{cor_definition}
Let $G$ be a compact Lie group and $T$ be a maximal torus.\ Then, the critical 
point set of any $T$-invariant Morse function on the flag manifold $G/T$ is 
given by its Weyl group. $\hspace{75mm}$ $\square$
\label{flag}\end{cor_definition}

\section{Intersections}\par

Let $G$ be a compact connected Lie group and $M$ be a compact smooth $G$-manifold. The following is our main result in this paper.
\begin{thm_definition}
Let $\Phi $ be a $G$-invariant Bott-Morse function on $M$. Let $p,q$ be $G$-fixed points. Assume the following conditions:\\
(1) $M^{G}\subset \operatorname{Cr}(\Phi )$.\\
(2) $\lambda (p)-\lambda (q)=2$.\\
(3) $W^{u}(p)$ and $W^{s}(q)$ intersect transversally.\\ 
Then every connected component of $\widetilde{\mathcal{M}}(p,q)$ is diffeomorphic to $S^{1}\times \mathbb{R}$.
\begin{proof}
 Let $C$ be a connected component of $\widetilde{\mathcal{M}}(p,q)$.\ By Lemma \ref{ggamma} and the connectedness of $G$, $C$ is a $G$-invariant subset of $\widetilde{\mathcal{M}}(p,q)$. We note that $C$ is non-compact. To see this, assume that $C$ is compact. Take $c'\in C$. Since the negative gradient flow $\gamma _{c'}(\mathbb{R})$ is connected, it must be contained in $C$. Therefore the assumption implies that $p=\lim\limits_{t\rightarrow -\infty}\gamma _{c}(t)\in C$.
This is a contradiction, because $p\not\in \widetilde{\mathcal{M}}(p,q)$. So $C$ is non-compact. Since $\operatorname{Cr}(\Phi )\cap C=\emptyset $, the assumption (1) implies that $M^{G}\cap C=\emptyset $. Let us show the following.\\
\\
\textbf{(3.1)} $\dim G\cdot c=1$.\\
\\
Assume that $\dim G\cdot c=2$. Then $G\cdot c$ is a codimension 0 submanifold of $C$. Therefore $G\cdot c$ is an open subset of $C$.
On the other hand,\ by the compactness of $G$,\ $G\cdot c$ is a closed subset of $C$. So we have $C=G\cdot c$ since $C$ is connected.
This is a contradiction, because $C$ is non-compact. Assume that $\dim G\cdot c=0$.\ Then by the connectedness of $G$, we have $G\cdot c=\{ c \}$.\ This is also a contradiction, because $c\not\in M^{G}$.\ Hence we 
have $\dim G\cdot c=1$. The proof of (3.1) is complete.

Define an action of $G\times \mathbb{R}$ on $C$ by $(g,t)\cdot c=g\cdot \gamma _{c}(t)$. In fact,
this gives an action on $C$, because
\begin{align*}
(gg',t+t')\cdot c & =gg'\cdot \gamma _{c}(t+t') \\  
                  & =g\cdot  \gamma _{g'\cdot c}(t+t') \\ 
                  & =g\cdot  \gamma _{\gamma _{g'\cdot c}(t')}(t) \\ 
                  & =(g,t)\cdot \gamma _{g'\cdot c}(t')  \\ 
                  & =(g,t)\cdot ((g',t')\cdot c) 
\end{align*}
for all $(g,t),(g',t')\in G\times \mathbb{R}$. We next show the following.\\
\\
\textbf{(3.2)} $(G\times \mathbb{R})_{c}=G_{c}\times \{ 0 \}$.\\
\\
Here, $(G\times \mathbb{R})_{c}$\ (resp. $G_{c}$) is the stabilizer of $c$ for the action of $G\times \mathbb{R}$ (resp. $G$) on $C$.
It is enough to show that $(G\times \mathbb{R})_{c}\subset G_{c}\times \{ 0\}$.\ Let $(g,t)$ be an element of $(G\times \mathbb{R})_{c}$.
It is sufficient to show $t=0$. Assume that $t>0$. Since $(g^{n},nt)\in (G\times \mathbb{R})_{c}$ for all $n\in \mathbb{N}$,
we have $\lim\limits_{n\rightarrow \infty}g^{n}\cdot c=\lim\limits_{n\rightarrow \infty} \gamma _{c}(-nt)=p$.  This implies that $p\in C$
since $G\cdot c$ is a closed subset of $C$. This is a contradiction.\ 
If we assume that $t<0$, a similar argument implies the same contradiction. 
The proof of (3.2) is complete.
\par
Let us consider the natural embedding $G\times \mathbb{R}/(G\times \mathbb{R})_{c}\longrightarrow \widetilde{\mathcal{M}}(p,q)$. By (3.1) and (3.2), we have $\dim G\times \mathbb{R}/(G\times \mathbb{R})_{c}=
\dim \widetilde{\mathcal{M}}(p,q)=2$. Thus $G\cdot \gamma _{c}(\mathbb{R})$ is open in $\widetilde{\mathcal{M}}(p,q)$.
In particular, every orbit of the action of $G\times \mathbb{R}$ on $C$ is open.
Since $C$ is connected, this implies that $C=G\cdot \gamma _{c}(\mathbb{R})$.
Therefore we obtain the following isomorphisms:
\begin{center}
$C\cong G\times \mathbb{R}/G_{c}\times \{0 \}\cong G/G_{c}\times \mathbb{R}\cong G\cdot c\times \mathbb{R}$.
\end{center}
By (3.1), $G\cdot c$ is a compact connected 1-dimensional manifold. Thus $G\cdot c$ is diffeomorphic to
$S^{1}$. Hence $C$ is diffeomorphic to $S^{1}\times \mathbb{R}$.

The proof is complete.
\end{proof}
\label{mainthm}
\end{thm_definition}

\begin{cor_definition}
Let $\Phi $ be a $G$-invariant Morse-Smale function on $M$. Let $p,q$ be critical points of $\Phi $ such that $\lambda (p)-\lambda (q)=2$. 
If $M^{G}$ is a finite set, every connected component of $\widetilde{\mathcal{M}}(p,q)$ is diffeomorphic to $S^{1}\times \mathbb{R}$.
\begin{proof}
By Proposition \ref{crit}, we have $M^{G}=\operatorname{Cr}(\Phi )$. So this corollary follows from Theorem \ref{mainthm}.
\end{proof}
\label{mainmorse}\end{cor_definition}
In the rest of this section, we study the stabilizer $G_{c}$. Let $G$ be a  compact connected Lie group which acts smoothly on $S^{1}$. We denote by $\mathfrak{g}$ the Lie algebra of $G$. Consider the following commutative diagram:
\begin{align*}
\xymatrix{ \ar@{}[rd]
G \ar[r]  & \operatorname{Diff}(S^{1}) \\
\mathfrak{g} \ar[u] \ar[r] & \operatorname{\Gamma}(TS^{1}). \ar[u]}
\end{align*}
Here, vertical arrows are exponential maps and horizontal arrows are induced by
the action of $G$ on $S^{1}$. Since $G$ is a compact connected Lie group, the exponential map $\mathfrak{g}\longrightarrow G$ is surjective.
Thus the image of $G\longrightarrow \operatorname{Diff}(S^{1})$ is completely determined by the image of
$\mathfrak{g}\longrightarrow \operatorname{\Gamma}(TS^{1})$. We need the following result of Plante [P, Theorem 1.2].

\begin{lem_definition}
Let $G$ be a Lie group and $\mathfrak{g}$ be its Lie algebra.
Assume that $G$ acts smoothly and transitively on $S^{1}$. Then the image of $\mathfrak{g}\longrightarrow \operatorname{\Gamma}(TS^{1}) $ 
is conjugate via a diffeomorphism to one of the following subalgebras of $\operatorname{\Gamma}(TS^{1}) $ \\
(1)$\displaystyle \left<\frac{\partial }{\partial x} \right>$,\\
(2)$\displaystyle \left< (1+\cos x)\frac{\partial }{\partial x},(\sin x)\frac{\partial }{\partial x}, (1-\cos x)\frac{\partial }{\partial x}\right>$. $\hspace{47mm}$
$\square$\\
\label{plante}\end{lem_definition}
Note that we have the isomorphism 
\begin{center}
$\displaystyle \left< (1+\cos x)\frac{\partial }{\partial x},(\sin x)\frac{\partial }{\partial x}, (1-\cos x)\frac{\partial }{\partial x}\right>\cong \mathfrak{sl}_{2}(\mathbb{R})$
\end{center}
of Lie algebras.
\begin{pro_definition}In the setting of Theorem \ref{mainthm}, let $C$ be a connected component of $\widetilde{\mathcal{M}}(p,q)$.
Then there is a surjective group homomorphism $\alpha :G\longrightarrow S^{1}$ and a diffeomorphism $C\cong S^{1}\times \mathbb{R}$ such that
the action of $G\times \mathbb{R}$ on $C\cong S^{1}\times \mathbb{R}$ is given by
\begin{center}
$(g,t)\cdot (x,s)=(\alpha (g)x,t+s)$
\end{center}
for all $(g,t)\in G\times \mathbb{R},(x,s)\in S^{1}\times \mathbb{R}$.
\begin{proof}
Take $c\in C$. We consider the action of $G$ on $G\cdot c$ and identify $G\cdot c$ with $S^{1}$.  Let $\alpha_{0} :G\longrightarrow \operatorname{Diff}(S^{1})$ be the representation of 
the action of $G$ on $S^{1}$, $\alpha_{0}': \mathfrak{g}\longrightarrow \operatorname{\Gamma}(TS^{1})$ the corresponding Lie algebra homomorphism.
  
 Since $\mathfrak{g}$ is the Lie algebra of the compact Lie group $G$, it does not admit $\mathfrak{sl}_{2}$ as a quotient Lie algebra . Hence by Lemma 
\ref{plante} we can take $\varphi  \in \operatorname{Diff}(S^{1})$
such that
\begin{center}
$\displaystyle \varphi _{*}(\alpha_{0}'(\mathfrak{g}))= \left<\frac{\partial }{\partial x} \right>$.
\end{center}

This shows that $\varphi (\alpha_{0}(G)) \varphi ^{-1} $  consists of rotations of $S^{1}$. 
Now we define a group homomorphism $\alpha :G\longrightarrow S^{1}$ by $\alpha (g):=\varphi\circ  \alpha_{0}(g)\circ \varphi ^{-1}$.
This map satisfies the required properties.
\end{proof}
\label{action}\end{pro_definition}
\begin{cor_definition}
In the setting of Theorem \ref{mainthm}, let $C$ be a connected component of $\widetilde{\mathcal{M}}(p,q)$. Then the stabilizer of $c\in C$ is independent of choice of $c$ and is a codimension 1 closed normal Lie subgroup of $G$.
$\hspace{42mm}$ $\square$
\end{cor_definition}

\end{document}